\documentclass[12pt]{article}

\usepackage{amsmath}
\usepackage{amsfonts}
\usepackage{amsthm}
\usepackage{graphicx}
\newtheorem{theorem}{Theorem}

\def\F{\mathcal{F}}
\def\C{\mathcal{C}}

\begin{document}

\title{Universal Cycles for Permutations}

\author{J Robert Johnson\\School of Mathematical Sciences\\Queen Mary,
  University of London\\Mile End Road, London E1 4NS, UK\\\emph{Email:} r.johnson@qmul.ac.uk}

\maketitle

\begin{abstract}
A universal cycle for permutations is a word of length $n!$ such that each of
the $n!$ possible relative orders of $n$ distinct integers occurs as a
cyclic interval of the word. We show how to construct such a universal
cycle in which only $n+1$ distinct integers are used. This is best
possible and proves a conjecture of Chung, Diaconis and Graham. 
\end{abstract}


\section{Introduction}

A \emph{de Bruijn cycle} of order $n$ is a word in $\{0,1\}^{2^n}$ in
which each $n$-tuple in $\{0,1\}^n$ appears exactly once as a cyclic
interval (see \cite{dB}).  The idea of a universal cycle generalizes the notion of a de Bruijn
cycle.

Suppose that $\F$ is a family of
combinatorial objects with $|\F|=N$, each of which is represented (not necessarily in a unique way) by an
$n$-tuple over some alphabet $A$. A \emph{universal cycle} (or \emph{ucycle})
for $\F$ is a word $u_1u_2\dots{u}_N$ with each
$F\in\F$ represented by exactly one $u_{i+1}u_{i+2}\dots{u}_{i+n}$
where, here and throughout, index addition is interpreted modulo $N$. With this terminology
a de Bruijn cycle is a ucycle for words of length $n$ over $\{0,1\}$
with a word represented by itself. The definition of ucycle was introduced by
Chung, Diaconis and Graham in \cite{CDG}. Their paper, and the
references therein, forms an good overview of the topic of universal cycles. The cases considered by them
include $\F$ being the set of permutations of an $n$-set,
$r$-subsets of an $n$-set, and partitions of an $n$-set. 

In this paper we will be concerned with ucycles for permutations: our
family $\F$ will be $S_n$, which we will regard as the set
of all $n$-tuples of distinct elements of $[n]=\{1,2,\dots{n}\}$. It is not immediately obvious how we should
represent permutations with words. The most natural thing to do would be to
take $A=[n]$ and represent a permutation by itself, but it is easily verified that
(except when $n\leq2$) it is not possible to have a ucycle in this
case. Indeed, if every cyclic interval of a word is to represent a
permutation then our word must repeat with period $n$, and so only $n$
distinct permutations can be represented. Another possibility which we
mention in passing would be
to represent the permutation $a_1a_2\dots{a}_n$ by
$a_1a_2\dots{a}_{n-1}$. It is clear that the permutation is determined
by this. It was shown by Jackson \cite{J} (using similar techniques to
those used for de Bruijn cycles) that these ucycles exist for all
$n$. Recently an efficient algorithm for constructing such ucycles
was given by Williams \cite{W}. He introduced the term \emph{shorthand
  universal cycles for permutations} to describe them. Alternatively, Chung, Diaconis and Graham in \cite{CDG}  
consider ucycles for permutations using a larger alphabet where each
permutation is represented by any $n$-tuple in which the elements have
the same relative order. Our aim is to prove their conjecture that such
ucycles always exist when the alphabet is of size $n+1$, the smallest
possible. In contrast to the situation with shorthand universal
cycles, the techniques used for de Bruijn cycles do not seen to help
with this so a different approach is needed.

To describe the problem more formally we need the notion of
order-isomorphism. If $a=a_1a_2\dots{a}_n$ and $b=b_1b_2\dots{b}_n$
are $n$-tuples of distinct integers, we say that $a$ and $b$
are $\emph{order-isomorphic}$ if 
\[
a_i<a_j \Leftrightarrow b_i<b_j
\]
for all $1\leq{i},j\leq{n}$. Note that no two distinct permutations in
$S_n$ are order-isomorphic, and that any $n$-tuple of distinct integers is order-isomorphic to exactly one permutation in
$S_n$. Hence, the set of $n$-tuples of distinct integers is
partitioned into $n!$ order-isomorphism classes which correspond to
the elements of $S_n$. 

We say that a word $u_1u_2\dots{u}_{n!}$ over an alphabet $A\subset\mathbb{Z}$ is a ucycle for $S_n$ if
there is exactly one $u_{i+1}u_{i+2}\dots{u}_{i+n}$ order-isomorphic
to each permutation in $S_n$. For example $012032$ is a ucycle for
$S_3$. Let $M(n)$ be the smallest integer $m$
for which there is a ucycle for $S_n$ with $|A|=m$. Note that if
$|A|=n$ then each permutation is represented by itself and so, as we
noted earlier, no ucycle is possible (unless $n\leq2$). We deduce that
$M(n)\geq{n+1}$ for all $n\geq3$. Chung, Diaconis and Graham in \cite{CDG} give the upper
bound $M(n)\leq{6n}$ and conjecture that $M(n)=n+1$ for all $n\geq3$. Our main result is that this
conjecture is true.

\begin{theorem}
For all $n\geq3$ there exists a word of length $n!$ over the alphabet
$\{0,1,2,\dots,n\}$ such that each element of $S_n$ is
order-isomorphic to exactly one of the $n!$ cyclic intervals of length $n$.
\end{theorem} 

We prove this by constructing such a word inductively. The details of
our construction are in the next section. Having shown that such a
word exists, it is natural to ask how many there are. In the
final section we give some bounds on this.

Our construction works for $n\geq5$. For smaller values of $n$ it is a
relatively simple matter to find such words by hand. For completeness
examples are $012032$ for $n=3$, and $012301423042103421302143$ for $n=4$.

\section{A Construction of a Universal Cycle}

We will show how to construct a word of length $n!$ over the alphabet
$\{0,1,2,\dots,n\}$ such that for each $a\in{S}_n$ there is a  cyclic
interval which is order-isomorphic to $a$. 

Before describing the construction we make a few preliminary definitions.

As is standard for universal cycle problems we let $G_n=(V,E)$, the
\emph{transition graph}, be the
directed graph with 
\begin{align*}
V&=\{(a_1a_2\dots{a_n}): a_i\in\{0,1,2,\dots,n\}, \text{ and } a_i\not=a_j
\text{ for all } i\not=j\}\\
E&=\{(a_1a_2\dots{a_n})(b_1b_2\dots{b_n}): a_{i+1}=b_i \text{ for all
} 1\leq{i}\leq{n-1}\}.
\end{align*}

Notice that every vertex of $G_n$ has out-degree and in-degree both
equal to 2.

The vertices on a  directed cycle in $G_n$ plainly correspond to the $n$-tuples which
occur as cyclic intervals of some word. Our task, therefore, is to find a directed cycle in $G_n$ of length $n!$ such that for each
$a\in{S_n}$ there is some vertex of our cycle which is
order-isomorphic to $a$. This is in contrast to many universal
cycle problems where we seek a Hamilton cycle in the transition graph.

We define the map on the integers:
\[
s_x(i)=\left\{
\begin{array}{ll}
i & \text{ if } i<x\\
i+1 &\text{ if } i\geq{x}.\\
\end{array}\right.
\]

We also, with a slight abuse of notation, write $s_x$ for the map
constructed by applying this map coordinatewise to an $n$-tuple. That is,
\[ 
s_x(a_1a_2\dots{a_n})=s_x(a_1)s_x(a_2)\dots{s_x}(a_n).
\]

The point of this definition is that if $a=a_1a_2\dots{a_n}\in{S}_n$
is a permutation of $[n]$ and $x\in[n+1]$ then $s_x(a)$ is the unique $n$-tuple of elements of
$[n+1]\setminus\{x\}$ which is order-isomorphic to $a$. Note that, as
will become clear, this is the definition we need even though our final construction will produce a ucycle for
permutations of $[n]$ using alphabet $\{0,1,2,\dots,n\}$.

We also define a map $r$ on $n$-tuples which permutes the elements of the
$n$-tuple cyclically. That is,
\[
r(a_1a_2\dots{a_n})=a_2a_3\dots{a}_{n-1}a_na_1.
\]

Note that $(a,r(a))$ is an edge of $G_n$ and that $r^n(a)=a$. 

As indicated above, we prove Theorem 1 by constructing a cycle of length $n!$ in $G_n$ such that for
each $a\in{S}_n$ the cycle contains a vertex which is order-isomorphic
to $a$. Our approach is to find a collection of short cycles in $G_n$
which between them contain one vertex from each order-isomorphism
class and to join them up. The joining up of the short cycles requires
a slightly involved induction step which is where the main work lies.

\begin{proof}[Proof of Theorem 1:]

\noindent
\emph{Step 1:} Finding short cycles in $G_n$

The first step is to find a collection of short cycles (each of length
$n$) in $G_n$ which between them contain exactly one element from each
order-isomorphism class of $S_n$. These cycles will use only $n$
elements from the alphabet and we will think of each cycle as being
``labelled'' with the remaining unused element. Suppose that for each
$a=a_1a_2\dots{a_{n-1}}\in{S}_{n-1}$ we choose a label $l(a)$ from
$[n]$. Let $0a$ be the $n$-tuple $0a_1a_2\dots{a_{n-1}}$. We have
the following cycle in $G_n$:
\[
s_{l(a)}(0a),r(s_{l(a)}(0a)),r^2(s_{l(a)}(0a)),\dots,r^{n-1}(s_{l(a)}(0a)).
\]

We denote this cycle by $\C(a,l(a))$. As an example, $\C(42135,2)$ is the
following cycle in $G_6$:
\[
053146\rightarrow531460\rightarrow314605\rightarrow146053\rightarrow460531\rightarrow605314\rightarrow053416,
\]
where arrows denote directed edges of $G_6$.

Note that for any choice of labels (that is any map $l$) the cycles
$\C(a,l(a))$ and $\C(b,l(b))$ are disjoint when
$a,b\in{S}_n$ are distinct. Consequently, whatever the choice of
labels, the collection of cycles
\[
\bigcup_{a\in{S}_{n-1}}\C(a,l(a))
\]
is a disjoint union. It is easy to see that the vertices on these
cycles contain between them exactly one $n$-tuple
order-isomorphic to each permutation in $S_n$.

We must now show how, given a suitable labelling, we can join up these short cycles.
\vskip5mm

\noindent
\emph{Step 2:} Joining two of these cycles

Suppose that $\C_1=\C(a,x)$ and
$\C_2=\C(b,y)$ are two of the cycles in $G_n$ described
above. What conditions on $a,b$ and their labels $x,y$ will allow us to join
these cycles?

We may assume that $x\leq y$. Suppose further that $1\leq x\leq{y}-2\leq n-1$, and that $a$ and $b$ satisfy the
following:
\[
b_i=\left\{
\begin{array}{ll}
a_i & \text{ if } 1\leq{a_i}\leq{x}-1\\
a_i+1 & \text{ if } x\leq{a_i}\leq{y}-2\\
x & \text{ if } a_i=y-1\\
a_i & \text{ if } y\leq{a_i}\leq{n-1}.\\
\end{array}\right.
\]

If this happens we will say that the pair of cycles $\C(a,x),\C(b,y)$
are \emph{linkable}.

In this case $s_{x}(0a)$ and $s_{y}(0b)$ agree at all but one
position; they differ only at the $t$ for which $a_t=y-1$ and
$b_t=x$. It follows that there is a directed edge in $G_n$ from
\[
r^t(s_{x}(0a))=s_x(a_t\dots{a}_{n-1}0a_1\dots{a}_{t-1})
\]
to
\[
r^{t+1}(s_{y}(0b))=s_y(b_{t+1}\dots{b}_{n-1}0b_1\dots{b}_t).
\]

Similarly, there is a directed edge in $G_n$ from
\[
r^t(s_{y}(0b))=s_x(b_t\dots{b}_{n-1}0b_1\dots{b}_{t-1})
\]
to
\[
r^{t+1}(s_{x}(0a))=s_y(a_{t+1}\dots{a}_{n-1}0a_1\dots{a}_t).
\]

If we add these edges to $\C_1\cup\C_2$ and remove
the edges
\[
r^t(s_{x}(0a))r^{t+1}(s_{x}(0a))
\]
and
\[
r^{t}(s_{y}(0b))r^{t+1}(s_{y}(0b)),
\]
then we produce a single cycle of length $2n$ whose vertices are
precisely the vertices in $\C_1\cup\C_2$. 

We remark that if $x=y-1$ then the other conditions imply that $a=b$
and so although we can perform a similar linking operation it is not
useful. If $x=y$ then $b$ is not well-defined.

As an example of the linking operation consider the linkable pair of 6-cycles $\C(42135,2)$ and
$\C(23145,5)$ in $G_6$. If we add the edges $531460\rightarrow314602$ and
$231460\rightarrow314605$, and remove the edges $531460\rightarrow314605$
and $231460\rightarrow314602$ then a single cycle of length 12 in $G_6$
is produced. These cycles and the linking operation are shown in
Figure 1.

\begin{figure}

\includegraphics[scale=0.75]{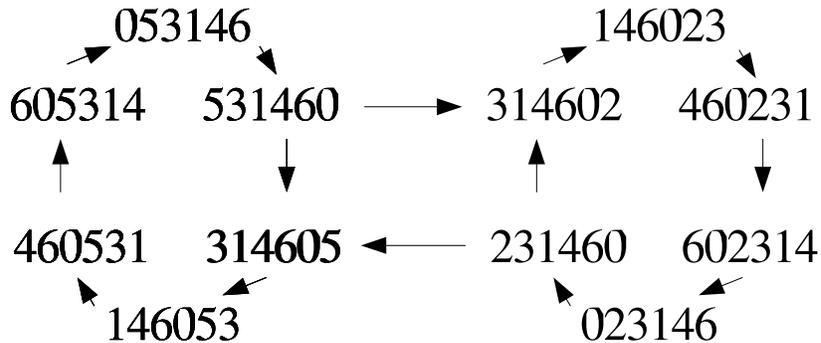}

\caption{The cycles $\C(42135,2)$ and $\C(23145,5)$ are linkable}

\end{figure}

\vskip5mm

\noindent
\emph{Step 3:} Joining all of these cycles

We now show that this linking operation can be used repeatedly to join
a collection of disjoint short cycles, one for each $a\in{S}_n$, together.

Let $H_{n}=(V,E)$ be the (undirected) graph with,
\begin{align*}
V&=\{(a,x):a\in{S}_{n-1}, x\in[n]\}\\
E&=\{(a,x)(b,y): \C(a,x), \C(b,y) \text{ are linkable}\}
\end{align*}

If we can find a subtree $T_{n}$ of $H_{n}$ of order $(n-1)!$ which contains
exactly one vertex $(a,x)$ for each $a\in{S}_{n-1}$ then we will be
able to construct the required cycle. Take any vertex $(a,l(a))$ of $T_{n}$ and consider the cycle
$\C(a,l(a))$ associated with it. Consider also the cycles associated
with all the neighbours in $T_n$ of $(a,l(a))$. The linking operation
described above can be used to join the cycles associated with these neighbours to $\C(a,l(a))$. This is
because the definition of adjacency in $H_{n}$ guarantees that we can join
each of these cycles individually. Also, the fact that every vertex in
$G_n$ has out-degree 2 means that the joining happens at different
places along the cycle. That is if $(b,l(b))$ and $(c,l(c))$ are
distinct neighbours of $(a,l(a))$ then the edge of $\C(a,l(a))$ which must
be deleted to join $\C(b,l(b))$ to it is not the same as the one which
 must be deleted to join $\C(c,l(c))$ to it. We conclude that we can
 join all of the relevant cycles to the cycle associated with $(a,l(a))$. The
connectivity of $T_n$ now implies that we can join all of the cycles
associated with vertices of $T_n$ into one cycle. This is plainly a
cycle with the required properties.

The next step is to find such a subtree in $H_{n}$.
\vskip5mm

\noindent
\emph{Step 4:} Constructing a Suitable Tree

We will prove, by induction on $n$, the stronger statement that for all $n\geq5$, there is a
subtree $T_n$ of $H_n$ of order $(n-1)!$ which satisfies:
\begin{enumerate}
\item for all $a\in{S}_{n-1}$ there exists a unique $x\in[n]$ such that $(a,x)\in{V}(T_n)$,
\item $(12\dots(n-1),1)\in{V}(T_n)$,
\item $(23\dots(k-1)1(k)(k+1)\dots(n-1),k)\in{V}(T_n) $ for all $3\leq{k}\leq{n}$,
\item $(32145\dots(n-1),2)\in{V}(T_n)$,
\item $(243156\dots(n-1),3)\in{V}(T_n)$,
\item $v(31245\dots(n-1))$ is a leaf in $T_n$,
\item $v(24135\dots(n-1))$ is a leaf in $T_n$.
\end{enumerate}
 
Where, for a tree satisfying property 1, we denote the unique vertex in $V(T_n)$ of the form $(a,x)$ by $v(a)$.

For $n=5$ a suitable tree can be found. One such is given in Figure 2.

\begin{figure}

\includegraphics[scale=0.75]{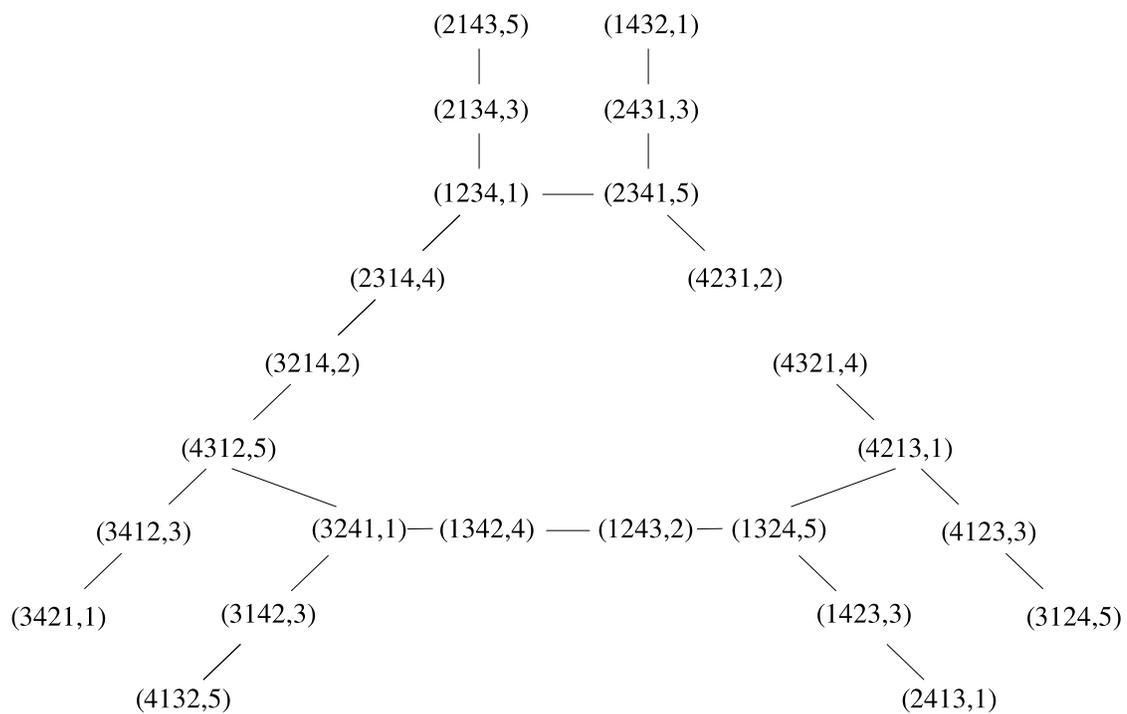}

\caption{A suitable choice for $T_5$}

\end{figure}

Suppose that $n\geq5$ and that we have a
subtree $T_n$ of the graph $H_n$ which satisfies the above
conditions. We will use this to build a suitable subtree of $H_{n+1}$.

A key observation for our construction is that the map from $V(H_n)$
to $V(H_{n+1})$ obtained by replacing each vertex
$(a_1a_2\dots{a}_{n-1},x)\in{V}(H_n)$ by
$(1(a_1+1)(a_2+1)\dots(a_{n-1}+1),x+1)\in{V}(H_{n+1})$ preserves
adjacency. It follows that subgraphs of $H_n$ are mapped into isomorphic
copies in $H_{n+1}$ by this map. Further, applying a fixed permutation to the coordinates of
the $n$-tuple associated with each vertex of $H_{n+1}$ gives an
automorphism of $H_{n+1}$ and so subgraphs of $H_{n+1}$ are mapped
into isomorphic copies. 

We take $n$ copies of $T_n$. These copies will be modified to form
the building blocks for our subtree of $H_{n+1}$ as follows. 

In the first copy we replace each vertex $(a_1a_2\dots{a_{n-1}},x)$ by
\[
(1(a_3+1)(a_1+1)(a_4+1)(a_2+1)(a_5+1)(a_6+1)\dots(a_{n-1}+1),x+1).
\]
By the observation above this gives an isomorphic copy of $T_n$ in $H_{n+1}$. We denote this tree by $T_{n+1}^{(0)}$.

In the next copy we replace each vertex
$(a_1a_2\dots{a_{n-1}},x)$ by
\[
((a_3+1)1(a_2+1)(a_1+1)(a_4+1)(a_5+1)\dots(a_{n-1}+1),x+1).
\]
We denote this tree by $T_{n+1}^{(1)}$.

For all $2\leq{k}\leq{n-1}$, we take a new copy of $T_n$ and modify it
as follows. We replace each vertex
$(a_1,a_2\dots{a}_{n-1},x)$ by
\[
((a_{k}+1)(a_1+1)(a_2+1)\dots(a_{k-1}+1)1(a_{k+1}+1)\dots(a_{n-1}+1),x+1).
\]
We denote these trees by $T_{n+1}^{(2)}, T_{n+1}^{(3)},\dots, T_{n+1}^{(n-1)}$.

As we mentioned this results in $n$ subtrees of $H_{n+1}$. They are clearly
disjoint because the position in which 1 appears in the first
coordinate of each vertex is distinct for distinct trees. Let $F_{n+1}$ be the $n$ component subforest of $H_{n+1}$ formed by taking the union
of the trees $T_{n+1}^{(k)}$ for $0\leq{k}\leq{n-1}$.

It is also easy to see that for every $a\in{S}_n$ there is a vertex in $F_{n+1}$ of the form
$(a,x)$ for some $x\in[n+1]$. It remains to show that the $n$ components
can be joined up to form a single tree with the required properties.

Notice that it is a consequence of the construction of the
$T_{n+1}^{(k)}$ that if $v(l_1l_2\dots{l}_{n-1})$ is a leaf in $T_n$
then the following vertices are all leaves in $F_{n+1}$:
\begin{itemize}
\item $v(1(l_3+1)(l_1+1)(l_4+1)(l_2+1)(l_5+1)(l_6+1)\dots(l_{n-1}+1))$
\item $v((l_3+1)1(l_2+1)(l_1+1)(l_4+1)(l_5+1)\dots(l_{n-1}+1))$
\item $v((l_{k}+1)(l_1+1)(l_2+1)\dots(l_{k-1}+1)1(l_{k+1}+1)\dots(l_{n-1}+1))$
  for $2\leq{k}\leq{n-1}$.
\end{itemize} 

We claim that $v(12\dots{n})$ is a leaf in $T_{n+1}^{(0)}$, and hence
a leaf in $F_{n+1}$. This follows from the fact that $v(24135\dots(n-1))$ is a leaf in $T_n$ and
the remark above. We delete this vertex from $F_{n+1}$ to form
a new forest. 

The construction of the $T_{n+1}^{(k-1)}$ and the fact that
$(23\dots(k-1)1(k)(k+1)\dots(n-1),k)\in{V}(T_n) $ for all
$3\leq{k}\leq{n}$ means that
$(23\dots{k}1(k+1)(k+2)\dots{n}),k+1)\in{V}(F_{n+1})$ for all
$3\leq{k}\leq{n}$. Further, the construction of the $T_{n+1}^{(1)}$
and the fact that $(32145\dots(n-1),2)\in{V}(T_n)$ means that
$(2134\dots{n},3)\in{V}(F_{n+1})$. 

We add to $F_{n+1}$ a new vertex $(123\dots{n},1)$ (this replaces the
deleted vertex from $T_{n+1}^{(0)}$) and edges from this
new vertex to $(23\dots{k}1(k+1)(k+2)\dots{n},k+1)$ for all
$2\leq{k}\leq{n}$. The previous observation shows that all of these
vertices are in $F_{n+1}$, and it is easy to check, using the
definition of linkability, that the
added edges are in $H_{n+1}$. This new forest has only two components.    

Similarly, the inductive hypothesis that $v(31245\dots(n-1))$ is
a leaf in $T_n$ gives that $v(342156\dots{n})$ is a leaf in
$F_{n+1}$ (using the case $k=3$ of the observation on
leaves). We delete this leaf from the forest, replace it with a new vertex
$(342156\dots{n},1)$, and add edges from this new vertex to
$(341256\dots{n},3)$ and $(143256\dots{n},4)$. The construction of
$T_{n+1}^{(2)}$ and the fact that $(32145\dots(n-1),2)\in{V}(T_n)$
ensures that the first of these vertices is in $F_{n+1}$. The
construction of $T_{n+1}^{(0)}$ and the fact that $(243156\dots(n-1),3)\in{V}(T_n)$
ensures that the second of these vertices is in our modified $T_{n+1}^{(0)}$. 

These modifications to $F_{n+1}$ produce a forest of one component --
that is a tree. Denote this tree by $T_{n+1}$. We will be done if we
can show that $T_{n+1}$ satisfies the properties demanded.

Plainly, $T_{n+1}$ contains exactly one vertex of the form $(a,x)$ for each
$a\in{S}_n$. By construction $(12\dots{n},1)$, and
$(23\dots{t}1(t+1)(t+2)\dots{n},t+1)$ are vertices of $T_{n+1}$ for
all $2\leq{t}\leq{n}$. Hence the first three properties are satisfied.

The construction of $T_{n+1}^{(2)}$ and the fact that
$(123\dots(n-1),1)\in{V}(T_n)$ ensures that $(32145\dots{n},2)$ is a
vertex of $T_{n+1}$. The construction of $T_{n+1}^{(3)}$ and the fact that
$(32145\dots(n-1),2)\in{V}(T_n)$ ensures that $(243156\dots{n},3)$ is a
vertex of $T_{n+1}$. Hence properties 4 and 5 are satisfied.

Finally, the construction of $T_{n+1}^{(1)}$ and the fact that
$v(31245\dots(n-1))$ is a leaf in $T_n$ ensures that
$v(31245\dots{n})$ is a leaf in $F_{n+1}$. The modifications which
$F_{n+1}$ undergoes do not change this and so it is a leaf in
$T_{n+1}$. The construction of $T_{n+1}^{(2)}$ and the fact that
$v(31245\dots(n-1))$ is a leaf in $T_n$ ensures that
$v(241356\dots{n})$ is a leaf in $F_{n+1}$. Again, the modifications to
$F_{n+1}$ do not change this and so this vertex is still a leaf in
$T_{n+1}$. Hence properties 6 and 7 are satisfied.

We conclude that the tree $T_{n+1}$ has the
required properties. This completes the construction.

\end{proof}

\section{Bounds on the Number of Universal Cycles}

Having constructed a ucycle for $S_n$ over the alphabet
$\{0,1,\dots,n\}$ it is natural to ask how
many such ucycles exist. We will regard words which differ only by a
cyclic permutation as the same so we normalize our universal
cycles by insisting that the first $n$ entries give a word
order-isomorphic to $12\dots{n}$. We denote by $U(n)$ the number of words of length $n!$ over the alphabet
$\{0,1,2,\dots,n\}$ which contain exactly one cyclic interval
order-isomorphic to each permutation in $S_n$ and for which the first
$n$ entries form a word which is order-isomorphic to
$12\dots{n}$. There is a natural upper bound which is essentially exponential
in $n!$ based on the
fact that if we are writing down the word one letter at a time we have 2
choices for each letter. We can also show that there is
enough choice in the construction of the previous section to prove a
lower bound which is exponential in $(n-1)!$. It is slightly
surprising that our construction gives a lower bound which is this large. However, the upper and
lower bounds are still far apart and we have no idea where the true answer
lies.

\begin{theorem}
\[
420^{\frac{(n-1)!}{24}}\leq{U}(n)\leq{(n+1)}2^{n!-n}
\]
\end{theorem}

\begin{proof}
Suppose we write down our universal cycle one letter at a time. We
must start by writing down a word of length $n$ which is
order-isomorphic to $12\dots{n}$; there are $n+1$ ways of doing
this. For each of the next $n!-n$ entries we must not choose any of the
previous $(n-1)$ entries (all of which are distinct) and so we have 2
choices for each entry. This gives the required upper bound.

Now for the lower bound. We will give a lower bound on the number of subtrees of $H_n$ which satisfies the
conditions of step 4 of the proof of Theorem 1. It can be checked that
if a universal cycle comes from a subtree of $H_n$ in the way
described then the tree is determined by the universal cycle. It
follows that the number of such trees is a lower bound for $U(n)$.

Notice that if we have $t_n$ such subtrees of $H_n$ then we have at least $t_n^n$ such
subtrees of $H_{n+1}$. This is because in our construction we took $n$
copies of $T_n$ to build $T_{n+1}$ from and each different set of choices
yields a different tree. We conclude that the number of subtrees of $H_n$
satisfying the conditions is at least 
\[
t_5^{5\times6\times\dots\times{n-1}}=t_5^{\frac{(n-1)!}{24}}.
\]

Finally, we bound $t_5$. We modify the given $T_5$ be adding edges from $(1432,1)$ to $(2143,5)$, from $(3421,1)$ to $(4132,5)$ and from
$(4231,2)$ to $(4321,4)$. This graph is such that that any of its
spanning trees satisfies the properties for our $T_5$. It
can be checked that this graph has 420 spanning trees. This gives
the lower bound.
\end{proof}

\section{Acknowledgements}

This work was inspired by the workshop on Generalizations of de Bruijn Cycles and Gray Codes held at Banff in December
2004. I thank the organisers and participants of the workshop for a
stimulating and enjoyable week.


\begin{thebibliography}{9}
\bibitem{CDG} Chung,~F., Diaconis,~P., Graham,~R. (1993) Universal
  cycles for combinatorial structures. \emph{Discrete Math.} {\bf 110}
  43--59.
\bibitem{dB} de Bruijn,~N.\,G. (1946) A combinatorial problem. \emph{Nederl. Akad. Wetensch., Proc.} {\bf 49} 758--764.
\bibitem{J} Jackson,~B.\,W. (1993) Universal cycles for $k$-subsets and $k$-permutations. \emph{Discrete Math.} {\bf 117} 141--150.
\bibitem{W} Williams,~A.\,M. (submitted 2007) Shorthand Universal Cycles for Permutations. \emph{ACM-SIAM Symposium on Discrete
  Algorithms 2008}.
\end{thebibliography}
\end{document}